\documentclass{amsart}

\usepackage{amssymb,tikz}
\usepackage{amsthm}
\usepackage{amsmath,a4wide}
\usepackage{amsbsy}
\usepackage[all]{xy}
\usepackage{bm}
\usepackage{hyperref}
\usepackage{tikz}
\usetikzlibrary{calc}
\usepackage{pgfplots}
\usepackage{caption}

\newtheorem{thm}{Theorem}

\DeclareMathOperator{\di}{\operatorname{div}}
\DeclareMathOperator{\curl}{\operatorname{curl}}

\begin{document}

\title{Fast escape in incompressible vector fields}

\author{Stefan Steinerberger}

\address{Department of Mathematics, Yale University, 10 Hillhouse Avenue, New Haven, 06511 CT, USA}
\email{stefan.steinerberger@yale.edu}

\begin{abstract} Swimmers caught in a rip current flowing away from the shore are advised to swim orthogonally to the current to escape it. We
describe a mathematical principle in a similar spirit.
More precisely, we consider flows $\gamma$ in the plane induced by incompressible vector fields $\textbf{v}:\mathbb{R}^2 \rightarrow
\mathbb{R}^2$ satisfying $  c_1 < \|v\| < c_2.$  The length $\ell$ a flow curve $\dot \gamma(t) = \textbf{v}(\gamma(t))$ until $\gamma$ leaves a disk of radius 1 
centered at the initial position can be as long as $\ell \sim c_2/c_1$. The same is true for the orthogonal flow $\textbf{v}^{\perp} = (-\textbf{v}_2, \textbf{v}_1)$. We show that a combination does strictly better: there always exists a curve flowing first along $\textbf{v}^{\perp}$ and then along $\textbf{v}$ which escapes the unit disk before reaching the length $ \sqrt{4\pi c_2 / c_1}$. Moreover, if the escape length
of $\textbf{v}$ is uniformly $\sim c_2/c_1$, then the escape length of $\textbf{v}^{\perp}$ is uniformly $\sim 1$ (allowing for a fast escape from the current). We also
prove an elementary quantitative Poincar\'{e}-Bendixson theorem that seems to be new.

\end{abstract}

\maketitle

\section{Introduction}
\subsection{Introduction}
We study an interesting property of solutions of autonomous systems in the plane when the underlying vector field $\textbf{v}:\mathbb{R}^2 \rightarrow \mathbb{R}^2$ is incompressible. Suppose you find yourself blindfolded in a the middle of the ocean and want to get to a point that is distance at least 1 from your starting point. Being blindfolded, simply swimming straight or using external landmarks for navigation is not possible: you can only use the local currents of the ocean. If swimming with the stream does take a very long time to get you to distance 1 from your starting point, then swimming \textit{orthogonally} to the stream is guaranteed to be very effective and get you out very quickly. 

\begin{figure}[h!]
\begin{minipage}[t]{0.49\columnwidth}%
\begin{center}
\includegraphics[width = 4.2cm]{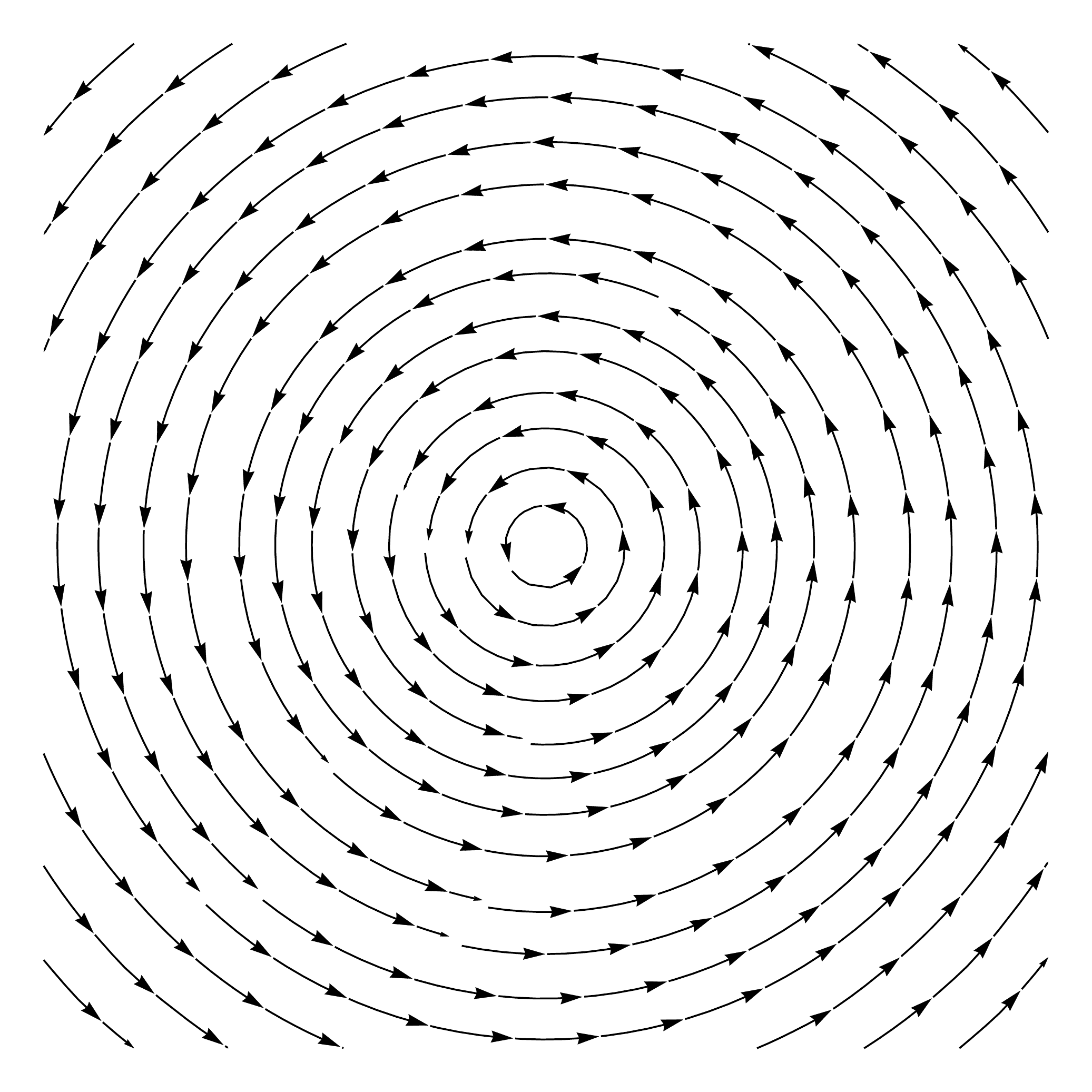}
\caption{Orbits of $(-y, x)$.}
\end{center}
\end{minipage}%
\begin{minipage}[t]{0.49\columnwidth}%
\begin{center}
\includegraphics[width = 4cm]{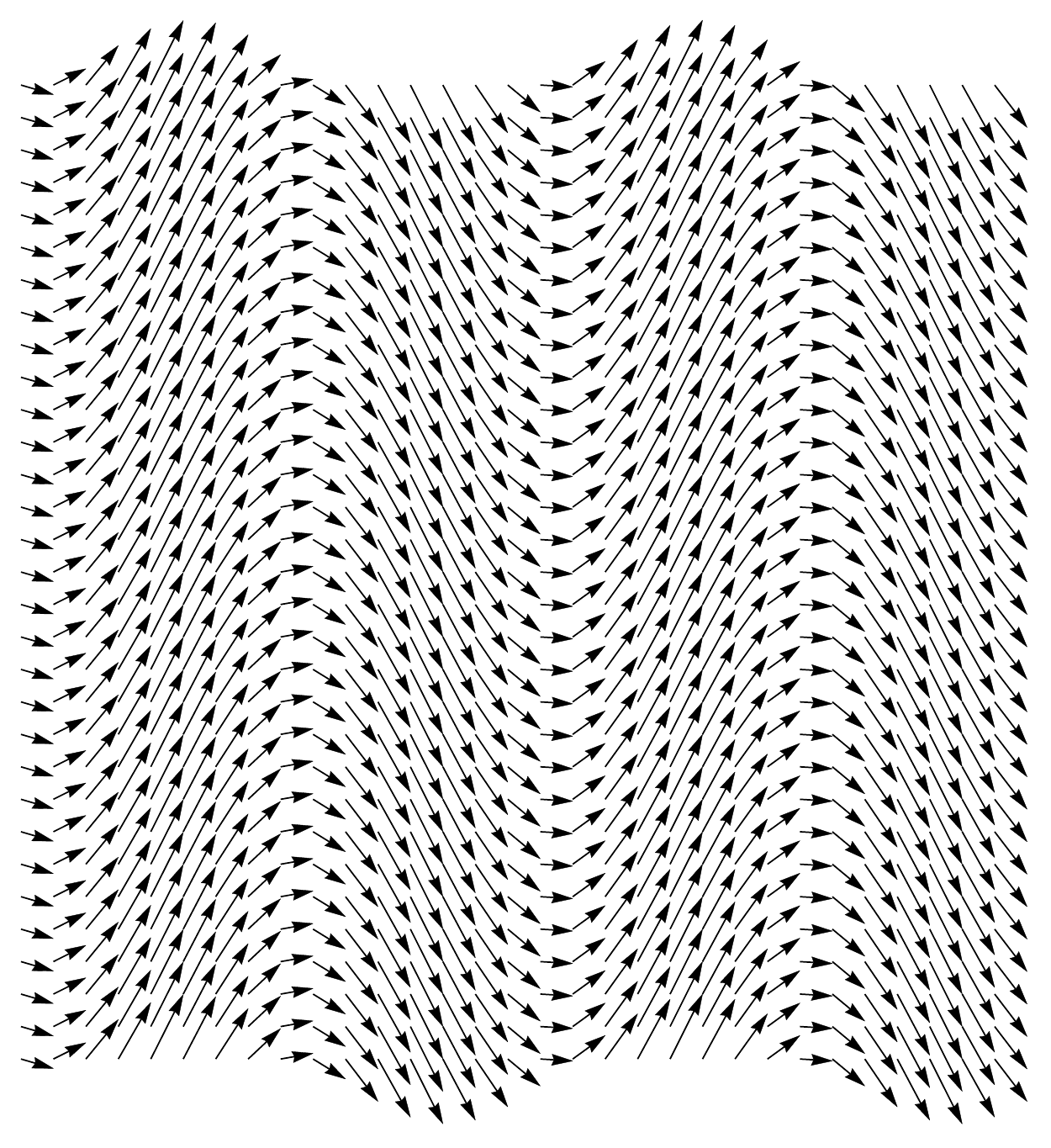}
\captionsetup{width=0.9\textwidth}
\caption{The example $\textbf{v}_N$ (see \S 3.1).}
\end{center}
\end{minipage}
\end{figure}

Put differently, a complicated global geometry of the stream lines of $\textbf{v} = (\textbf{v}_1, \textbf{v}_2)$ forces the flow along the orthogonal vector field $\textbf{v}^{\perp}=(-\textbf{v}_2, \textbf{v}_1)$ to be rather simple.  Needless to say, such a statement fails dramatically for general vector fields. More precisely, suppose we start in $\gamma(0) \in \mathbb{R}^2$ and are exposed to a smooth, incompressible vector field $\textbf{v}: \mathbb{R}^2 \rightarrow \mathbb{R}^2$, i.e.
$$ \frac{\partial}{\partial x}\textbf{v}_1 +  \frac{\partial}{\partial y}\textbf{v}_2 = 0.$$
It is entirely possible that following the vector field will lead us in circles, an easy example being 
$ \textbf{v}(x,y) = (-y, x).$
This is illustrated in Fig.1 where following the vector field will never allow an escape from the unit disk even if one starts outside of the origin (in which the vector field vanishes): to prevent this, we additionally assume throughout the paper that
$$ 0 < c_1 < \|\textbf{v}\| < c_2 < \infty.$$
The Poincar\'{e}-Bendixson theorem (see \cite{ben, poi, book}) then implies that the flow will eventually escape the unit disk:  if the flow curve of such a vector field were contained in a compact subset of $\mathbb{R}^2$, then 
it is either periodic or has to approach a limit cycle but any limit cycle encloses at least one critical point \cite[Corollary 1.8.5.]{guck} the existence of which would violate the assumption $\|\textbf{v}\| > c_1$. Our statements will be about the length of the flow curves as curves in the plane and not about the value of time $t$ in $\gamma(t)$. Such geometric quantities are invariant under dilating the vector field with a constant scalar
$$  (\textbf{v}_1, \textbf{v}_2) \rightarrow (\lambda \textbf{v}_1, \lambda \textbf{v}_2) \qquad \mbox{for some}~\lambda > 0$$
which corresponds to rescaling time $t \rightarrow t/\lambda$. Invariance under this scaling implies that our results about the geometry of curves will only depend on the ratio $c_2/c_1$.

\section{Main results}
Explicit examples, which are given in the next section, show that the length of the flow induced by either $\textbf{v}$ and $\textbf{v}^{\perp}$ 
might require a total length $\ell \sim c_2/c_1$ to escape the unit disk, where
$$ \ell = \inf_{t > 0}{ \left\{ \int_{0}^{t}{\|\gamma'(s)\|ds}: \|\gamma(t) - \gamma(0)\| = 1\right\}}.$$
We assume $\textbf{v}:\mathbb{R}^2 \rightarrow \mathbb{R}^2$ is a smooth, divergence-free vector field
satisfying $ c_1 < \| \textbf{v} \| < c_2$
for some $0 < c_1 < c_2 < \infty$. $\textbf{v}^{\perp}$ denotes the vector field $(-\textbf{v}_2, \textbf{v}_1)$. We write $e^{t\textbf{v}}$ to denote the 
flow induced by the vector field $\textbf{v}$ and $e^{t \textbf{v}^{\perp}}$ for the flow induced by $\textbf{v}^{\perp}$. We abuse notation and write the length of the flow line as
$$  |e^{t \textbf{v}^{}}\circ e^{s \textbf{v}^{\perp}} \gamma(0)| := \int_{0}^{s}{\left\| \frac{\partial}{\partial u}e^{u \textbf{v}^{\perp}} \gamma(0)\right\| du} + \int_{0}^{t}{ \left\| \frac{\partial}{\partial u}e^{u \textbf{v}} \left( e^{s \textbf{v}^{\perp}} \gamma(0) \right)\right\| du}$$
which is the length of curve given in $\mathbb{R}^2$ by concatenation of the flows
$$\gamma(0) \underbrace{\longrightarrow}_{\textbf{v}^{\perp}} e^{s \textbf{v}^{\perp}} \gamma(0)  \underbrace{\longrightarrow}_{\textbf{v}^{}}   e^{t \textbf{v}} \circ e^{s \textbf{v}^{\perp}} \gamma(0).$$

\begin{figure}[h!]
\begin{center}
\begin{tikzpicture}[scale = 1]
\draw[ thick] (0,0) to[out=30,in=180] (1,0.5)   to[out=0,in=180] (2,-0.5)     to[out=0,in=180] (3,0.5) to[out=0,in=180] (4,-0.5) to[out=0,in=180] (5,0.5) to[out=0,in=180] (6,-0.5);  
\draw[ thick] (0,0.5) to[out=30,in=180] (1,1)   to[out=0,in=180] (2,0)     to[out=0,in=180] (3,1) to[out=0,in=180] (4,0) to[out=0,in=180] (5,1) to[out=0,in=180] (6,0);  
\filldraw (0.5,0.3) circle (2pt);
\draw[dashed, thick] (0.9,-0.5) to[out=100,in=290] (-0.2,1.4);
\filldraw (0,1) circle (2pt);
\node at (-0.7,0.8) {$e^{s\textbf{v}^{\perp}}\gamma(0)$};
\filldraw (5,1.6) circle (2pt);
\node at (5.3,1.9) {$e^{t\textbf{v}^{}}e^{s\textbf{v}^{\perp}}\gamma(0)$};
\node at (0.4,-0.3) {$\gamma(0)$};
\draw[ thick] (0,1) to[out=10,in=170] (6,1.5); 
 \draw [dashed] (5.1,-0.6) arc (-10:18:5);
\end{tikzpicture}
\end{center}
\caption{The statement illustrated: starting in $\gamma(0)$, the flow curves of $\textbf{v}$ (thick) are long; flowing along $\textbf{v}^{\perp}$ (dashed) it is possible to reach a
short flow line of $\textbf{v}$.}
\end{figure}
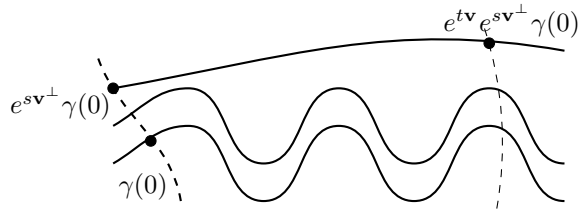

\begin{thm} Let $\textbf{v}:\mathbb{R}^2 \rightarrow \mathbb{R}^2$ be incompressible and smooth in a neighborhood of the disk of radius 1 around $\gamma(0) \in \mathbb{R}^2$. There exist $s,t \geq 0$ such that $\| (e^{t \textbf{v}} \circ e^{s \textbf{v}^{\perp}})\gamma(0) - \gamma(0)\| = 1$ and 
$$ |e^{t \textbf{v}} \circ e^{s \textbf{v}^{\perp}}\gamma(0)| \leq \sqrt{4\pi} \sqrt{ \frac{c_2}{c_1} }.$$
\end{thm}
If the conditions $s,t \geq 0$ are relaxed to $s,t \in  \mathbb{R}$ (which allows for both swimming with or against the current), then the constant can be improved by a factor of 2. 
There is an aspect to the Theorem not covered by Figure 5. It is certainly possible that \textit{all} flow lines of $\textbf{v}$ are equally inefficient
(see the next section for an example): in that case, the regularizing flow $\textbf{v}^{\perp}$ must necessarily be efficient itself
because the statement allows for $s=0$ or $t = 0$ (obviously not both at the same time). 

\begin{thm} Let $\textbf{v}:2 \mathbb{D} \rightarrow \mathbb{R}^2$ be incompressible and  $c_1 \leq \|\textbf{v}\| \leq c_2$. If, for all $\gamma(0) \in \mathbb{D}$,
$$  \inf_{t \in \mathbb{R}}{\left\{  | e^{t \textbf{v}}\gamma(0) |: \|e^{t \textbf{v}}\gamma(0) - \gamma(0)\| = 1\right\}} \geq \mathcal{L},$$
then
$$\inf_{t \in \mathbb{R}}{\left\{  | e^{t \textbf{v}^{\perp}}(0,0) |: \|e^{t \textbf{v}^{\perp}}(0,0)\| = 1\right\}} \leq \pi \frac{c_2 / c_1}{\mathcal{L}}.$$
\end{thm}
It is not a difficult to modify the proof to obtain a version of the statement, where we restrict $t > 0$ in both the assumption and the conclusion at the cost of a larger constant.\\

 We are not aware of any directly related results. However, any lifeguard knows about rip currents: these are strong outward flows that are highly localized but powerful enough to carry swimmers away from the shore. The common advice is to swim orthogonally to the flow to escape the (usually quite narrow) region
where they occur. 

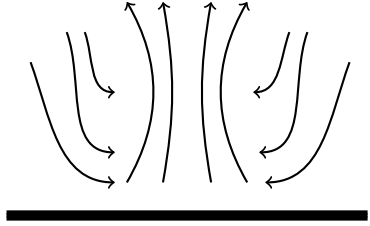
\begin{figure}[h!]
\centering
\begin{tikzpicture}[scale=0.8]
\draw [ultra thick] (0,0) -- (6,0);
\draw [ultra thick] (0,-0.1) -- (6,-0.1);
\draw [ultra thick] (0,-0.05) -- (6,-0.05);
\draw[thick, ->] (2,0.5) to[out=60,in=300] (2,3.5);
\draw[thick, ->] (2.6,0.5) to[out=80,in=280] (2.6,3.5);
\draw[thick, ->] (4,0.5) to[out=120,in=240] (4,3.5);
\draw[thick, ->] (3.4,0.5) to[out=100,in=260] (3.4,3.5);
\draw[thick, ->] (1.3,3) to[out=290,in=180] (1.8,2);
\draw[thick, ->] (1,3) to[out=290,in=180] (1.8,1);
\draw[thick, ->] (0.4,2.5) to[out=290,in=180] (1.8,0.5);
\draw[thick, ->] (5.7,2.5) to[out=250,in=0] (4.3,0.5);
\draw[thick, ->] (4.7,3) to[out=250,in=0] (4.1,2);
\draw[thick, ->] (5,3) to[out=250,in=0] (4.2,1);
\end{tikzpicture}
\caption{A rip current close to a beach}
\end{figure}

There is also an old problem of Bellman \cite{bell1, bell2, bell3} (see also the survey of Finch \& Wetzel \cite{finch})
\begin{quote}
A hiker is lost in a forest whose shape and dimensions are precisely known to him. What is the
best path for him to follow to escape from the forest?
\end{quote}
'Best' can be understood as minimizing the worst case length or the expected length. Our problem is, in a vague philosophical sense, dual: in Bellman's problem, the location is unknown but movement is unrestricted while in our problem the location is perfectly understood (at the center of a disk of radius 1) while movement is restricted and the way we can move at certain points is unknown (but, since its generated by an incompressible vectorfield, it does have some regularity). \\

Clearly, this result suggests a large number of variations --
an elementary variation we have been unable to find in the literature is the following quantitative Poincar\'{e}-Bendixson theorem. 

\begin{thm}[Quantitative Poincar\'{e}-Bendixson] Let $\textbf{v}:2\mathbb{D} \rightarrow \mathbb{R}^2$ be a smooth vector field satisfying
$c_1 < \|v\| < c_2$ and let $\gamma(t)$ be the flow induced by $\textbf{v}$ acting on $\gamma(0) = (0,0)$. Then 
$$ \inf_{t > 0}{ \left\{ \int_{0}^{t}{\|\gamma'(s)\|ds}: \|\gamma(t)\| = 1\right\}} \leq  \pi \frac{c_2}{c_1}+ \frac{1}{ c_1}  \int_{\mathbb{D}}{ \left| \curl\textbf{v}\right|~ dx dy} .$$
\end{thm}
The proof is not difficult
and similar in spirit to the proof of the classical Dulac criterion showing the nonexistence of closed orbits.
The statement captures the intuitive notion that in order for the flow to spend a long time in the domain, one requires a large amount
of rotation in the vector field. There is nothing special about the unit disk, for general simply-connected domains $\Omega$, the same result holds true
with $\pi$ replaced by $|\partial \Omega|/2$.
Note that the integral term vanishes completely for irrotational vector fields. In that case, there is a very easy proof: an irrotational $\textbf{v}$ can be written as a gradient
flow $\textbf{v} = \nabla A$. Since $\|  \textbf{v} \| = \|\nabla A\| \leq c_2$, 
$$ \sup_{x \in \mathbb{D}} A(x) - A(0,0) \leq c_2.$$
Using that $\textbf{v}$ generates the gradient flow of $A$, we may use that $A$, along that flow, increases at least as quickly as $c_1$ times the arclength measure. Hence,
$$   c_1 \inf_{t > 0}{ \left\{ \int_{0}^{t}{\|\gamma'(s)\|ds}: \|\gamma(t)\| = 1\right\}} \leq  c_2$$
which is stronger than Theorem 3 by a factor of $\pi$.\\

\textit{Open problems.} Even in two dimensions many natural questions remain. It is not clear to us how one would construct extremal configurations giving sharp constants. Theorem 1 guarantees the existence of an efficient route which first flows along $\textbf{v}^{\perp}$ and then along $\textbf{v}$. Suppose we allow for changing the vector field twice (i.e. first $\textbf{v}^{\perp}$, then $\textbf{v}$ and then once more $\textbf{v}^{\perp}$): can one reduce the length even further? How does the minimal length scale if one is allowed to switch $k-1$ times?
Similar questions can be asked in higher dimensions but the concept of orthogonal flow is no longer
uniquely defined -- it is not clear to us whether and in which form any of the results could be generalized and we believe this could be of some interest.

\section{Explicit constructions}

\subsection{Incompressible flows.} Consider (see Fig. 2) the incompressible vector field
$$ \textbf{v}_N(x,y) = \left(1, \frac{N}{2} \cos{\left(N x\right)}\right) \qquad \mbox{satisfying} \qquad 1 \leq \|\textbf{v}\| \leq \sqrt{\frac{N^2}{4}+1}.$$
Note that $c_2/c_1 \sim N$.
The initial condition $\gamma(0) = (0,0)$ gives rise to the flow $\gamma: \mathbb{R} \rightarrow \mathbb{R}^2$ 
$$ \gamma(t) = \left(t, \frac{1}{2}\sin{(Nt)}\right).$$
The distance of $\gamma(t)$ to the origin is bounded by
$$ \| \gamma(t) \| = \sqrt{t^2 + \frac{1}{4}\sin^2{(N t)}} \leq \sqrt{t^2 + \frac{1}{4}}.$$  
Therefore, in oder to have $\|\gamma(t)\| \geq 1$, we certainly need $t \geq \sqrt{3}/2$. The length of the curve up to that point is at least
$$ \int_{0}^{ \sqrt{3}/2}{\| \gamma'(t)\| dt} = \int_{0}^{ \sqrt{3}/2}{\sqrt{1 + \frac{N^2}{4}\cos^2{(Nt)}}dt} \geq  \frac{N}{2}\int_{0}^{ \sqrt{3}/2}{|\cos{(Nt)}|dt} \geq \frac{N}{8}.$$
This means that $\gamma(t)$ travels a total length of $\ell \sim N \sim c_2/c_1$ before leaving a disk of radius 1. 

\subsection{Irrotational flows.} The next natural question is whether it might be advantageous to immediately follow the orthogonal flow $\textbf{v}^{\perp} = (-\textbf{v}_2, \textbf{v}_1)$. It turns out that here, too, the
curve can be as long as $\ell \sim c_2/c_1$: in two
dimensions, we can use 
$$ \curl \textbf{v}^{\perp} = \curl  (-\textbf{v}_2, \textbf{v}_1) = \frac{\partial \textbf{v}_1}{\partial x} + \frac{\partial \textbf{v}_2}{\partial y} =  \di \textbf{v} = 0.$$
Hence, $\textbf{v}^{\perp}$ is irrotational and can therefore be written as $\textbf{v}^{\perp} = \nabla A$ for some scalar potential $A: \mathbb{R}^2 \rightarrow \mathbb{R}$. 
It suffices to construct an $A$ with $c_1 \leq \| \nabla A\| \leq c_2$ such that its gradient flow moves in a zigzag line. We construct
the scalar potential $A$ explicitely: given $N \in \mathbb{N}$  we define the curve $\gamma:\mathbb{R} \rightarrow \mathbb{R}^2$ via
$$ \gamma(t) = \left(\sin{(Nt)}, t\right).$$
Our (merely Lipschitz-continuous) candidate for $A$ is
$$ \tilde A(x,y) = Ny - 2N\inf_{t \in \mathbb{R}}{\| (x,y) - \gamma(t)\|}.$$
\begin{figure}[h!]
\begin{minipage}[t]{0.49\columnwidth}%
\begin{center}
\begin{tikzpicture}[scale = 1.6]
\draw[ultra thick] (0,0) to[out=10,in=270] (1,0.5)   to[out=90,in=350] (0,1)     to[out=170,in=270] (-1,1.5) to[out=90,in=190] (0,2);  
\node at (0.4,-0.3) {\Large $\gamma$};
\draw [ ->] (-0.1,0.3)--(0,0.1);
\draw [ ->] (0.4,0.5)--(0.7,0.5);
\draw [ ->] (-0.1,0.7)--(-0.05,0.9);

\draw [ ->] (0.05,1.3)--(0,1.1);
\draw [ ->] (-0.4,1.5)--(-0.7,1.5);
\draw [ ->] (0.1,1.7)--(0.05,1.9);

\draw [ ->] (-0.9,0.7)--(-0.8,1);
\draw [ ->] (0.9,1.3)--(0.8,1);
\end{tikzpicture}
\caption{A part of the curve $\gamma$ and the local gradients of $\tilde A$.}
\end{center}
\end{minipage}%
\begin{minipage}[t]{0.49\columnwidth}%
\begin{center}
\includegraphics[width = 7.2cm]{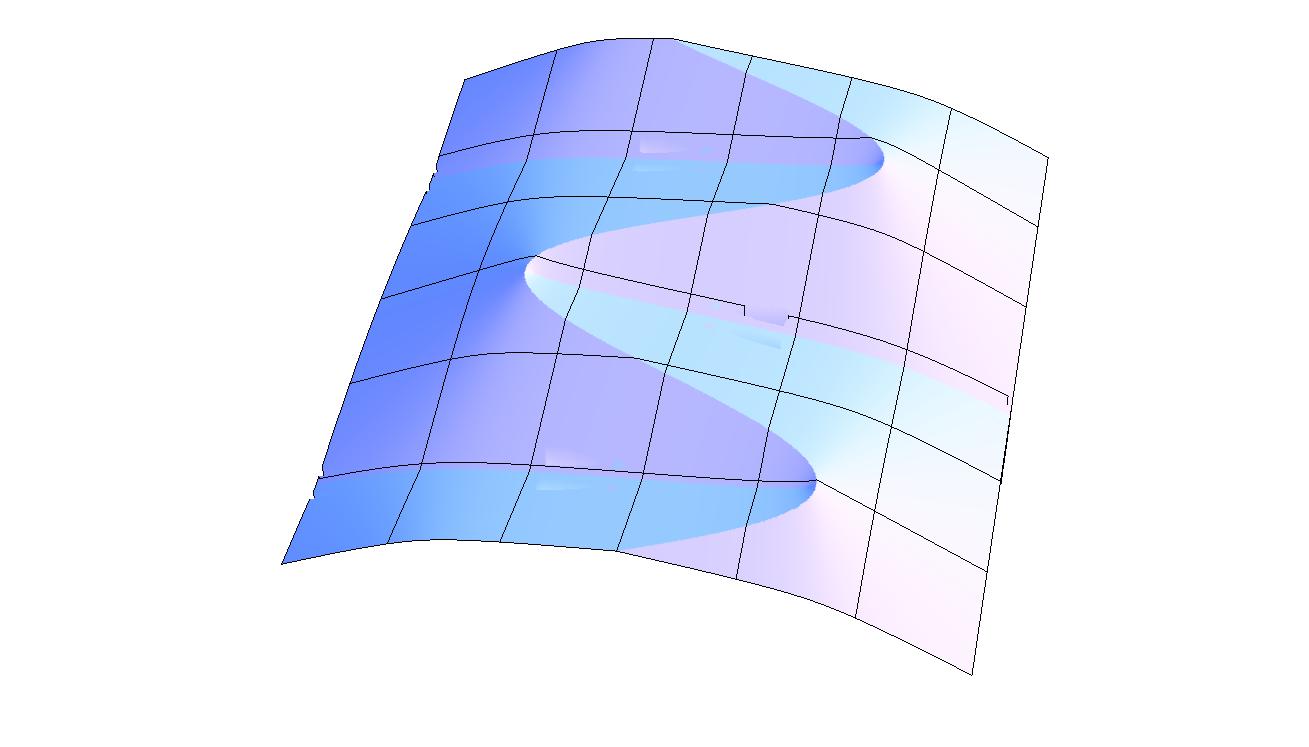}
\caption{A plot of $\tilde A$.}
\end{center}
\end{minipage}
\end{figure}
We remark that following $\gamma$ for $0 \leq t \leq 1$ gives rise to a curve of length $\sim N$. It is not difficult to see that, whenever defined, we have $0.1 \leq \|\nabla \tilde A \| \leq 10N$
and that, by the usual density arguments,  $\tilde A$ can be approximated in the $C^1-$norm by a smooth function $A$ having all necessary properties. This yields
an irrotational vector field $\textbf{v}^{\perp} = \nabla A$ where the escape length is of order $N \sim c_2/c_1$.

\section{Proofs}
\subsection{Proof of the Theorem 1.} The proof combines basic properties of incompressible vector fields in the plane with the coarea formula. In two dimensions we have 
$$ \di \textbf{v} = \curl \textbf{v}^{\perp} \quad \mbox{and thus} \quad \curl \textbf{v}^{\perp} = 0 \quad \mbox{for incompressible}~\textbf{v}.$$
Therefore there exists a scalar function $A:\mathbb{R}^2 \rightarrow \mathbb{R}$
such that
$$ \textbf{v}^{\perp} = \nabla A.$$
Trivially,
 $$\left\langle \textbf{v}, \nabla A \right\rangle = \left\langle \textbf{v}, \textbf{v}^{\perp} \right\rangle = 0.$$
Altogether, we have the existence of a function $A: \mathbb{R}^2 \rightarrow \mathbb{R}$ such that $\textbf{v}^{\perp}$ is
the gradient flow with respect to $A$ while $\textbf{v}$ flows along the level sets of the function $A$. We know that the function $A$ cannot change
too much on a set of large measure inside the unit disk $\mathbb{D}$ because
$$ \int_{\mathbb{D}}{| \nabla A |}  =  \int_{\mathbb{D}}{\| \textbf{v}^{\perp} \|} =   \int_{\mathbb{D}}{\| \textbf{v} \|}  \leq c_2 \pi.$$
We may assume w.l.o.g. (after possibly adding a constant) that $A(\gamma(0)) = 0$. The next ingredient is the coarea formula (see e.g. Federer \cite{federer})
$$  \int_{-\infty}^{\infty}{\mathcal{H}^1(A^{-1}(t))dt} =  \int_{\mathbb{D}}{| \nabla A |},$$
where $\mathcal{H}^1$ is the one-dimensional Hausdorff measure. Hence, for every $z > 0$, we have 
$$ z \inf_{0 < t < z}{\mathcal{H}^1(A^{-1}(t))} \leq \int_{0}^{z}{\mathcal{H}^1(A^{-1}(t))dt}  \leq \int_{-\infty}^{\infty}{\mathcal{H}^1(A^{-1}(t))dt} = \int_{\mathbb{D}}{| \nabla A |}  \leq c_2 \pi.$$
Setting $z = \sqrt{ \pi c_1 c_2}$, this implies the existence of a $0 < t_0 \leq \sqrt{c_1 c_2 \pi}$ with
$$ \mathcal{H}^1(A^{-1}(t_0)) \leq \sqrt{\frac{c_2}{c_1} \pi }.$$
\begin{figure}[h!]
\begin{center}
\begin{tikzpicture}[scale = 1]
\draw[ thick] (0,0) to[out=30,in=180] (1,0.5)   to[out=0,in=180] (2,-0.5)     to[out=0,in=180] (3,0.5) to[out=0,in=180] (4,-0.5) to[out=0,in=180] (5,0.5) to[out=0,in=180] (6,-0.5);  
\filldraw (0.5,0.3) circle (2pt);
\node at (8,1) {$\left\{\textbf{x}: A(\textbf{x}) = t_0 \right\}$};
\draw [thick, ->] (6.6, 1)  to[out=200,in=270]  (5.5,1.5);
\node at (8,0) {$\left\{\textbf{x}: A(\textbf{x}) = 0 \right\}$};
\draw [thick, ->] (6.6, 0)  to[out=160,in=20]  (5.7,0);
\draw[dashed, thick] (0.9,-0.5) to[out=100,in=290] (-0.3,1.5);
\filldraw (0,1) circle (2pt);
\node at (0.4,-0.3) {$\gamma(0)$};
\draw[ thick] (0,1) to[out=10,in=170] (6,1.5); 
\end{tikzpicture}
\end{center}
\caption{There exists a flow line of length at most $\sqrt{\pi c_2/c_1}$ described as level set $A(\textbf{x}) = t_0$ for some $|t_0| \leq \sqrt{\pi c_1c_2}$. It remains to estimate the length of the dashed line.}
\end{figure}
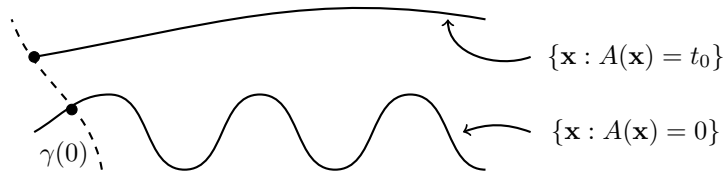

Let $s_0 > 0$ be the smallest number for which
$$ A\left( e^{s_0 \textbf{v}^{\perp} \gamma(0)} \right) = t_0.$$
A priori it might not clear that such a $s_0$ exists: it follows from the fact that $\textbf{v}^{\perp}$ is the gradient flow for $A$ and $\|\textbf{v}^{\perp}\| \geq c_1$. More precisely, using the fundamental theorem of calculus and the chain rule, we have
$$A\left( e^{s_0 \textbf{v}^{\perp}} \gamma(0) \right) = \int_{0}^{s_0}{\left(  A\left( e^{s \textbf{v}^{\perp}}\gamma(0) \right) \right)' ds} = \int_{0}^{s_0}{ \left\langle ( \nabla A)( e^{s \textbf{v}^{\perp}}\gamma(0)), (e^{s \textbf{v}^{\perp}}\gamma(0))' \right\rangle ds} .$$
However, since $\textbf{v}^{\perp}$ is the gradient flow of $A$, the scalar product simplifies for every $s$ to
$$  \left\langle  \nabla A \left( e^{s \textbf{v}^{\perp}} \gamma(0)\right), (e^{s \textbf{v}^{\perp}}\gamma(0))' \right\rangle = \left\| \nabla A \left( e^{s \textbf{v}^{\perp}} \gamma(0)\right)\right\|  \|(e^{s \textbf{v}^{\perp}}\gamma(0))'\|  \geq 
c_1 \|(e^{s \textbf{v}^{\perp}}\gamma(0))'\| $$
and therefore
$$ A\left( e^{s_0 \textbf{v}^{\perp}}\gamma(0) \right) =  \int_{0}^{s_0}{ \left\langle ( \nabla A)( e^{s \textbf{v}^{\perp}}\gamma(0)), (e^{s \textbf{v}^{\perp}}\gamma(0))' \right\rangle ds} \geq c_1  \int_{0}^{s_0}{ \|(e^{s \textbf{v}^{\perp}}\gamma(0))'\| ds} = c_1 | e^{s_0 \textbf{v}^{\perp}}\gamma(0)|.$$
This implies that the length from $\gamma(0)$ along $\textbf{v}^{\perp}$ until it hits $A^{-1}(t_0)$ is at most 
$$|e^{s_0 \textbf{v}^{\perp}}\gamma(0)| \leq \frac{A\left( e^{s_0 \textbf{v}^{\perp}}\gamma(0) \right) }{c_1} = \frac{t_0}{c_1} \leq \frac{\sqrt{c_1 c_2 \pi}}{c_1} = \sqrt{\pi}\sqrt{\frac{c_2}{c_1}}.$$
Summarizing, the curve from $(0,0)$ along $\textbf{v}^{\perp}$ until it intersects $A^{-1}(t_0)$ has length at most $\sqrt{c_2\pi/c_1}$ and the level set $\left\{ \textbf{x}: A(\textbf{x}) = t_0\right\}$ has length at most $\sqrt{c_2\pi/c_1}$ and this yields the statement. $\qed$

\subsection{Proof of Theorem 2.} Theorem 2 can be shown by a variant of the argument above. Let us first consider the curve
$$ \gamma = \left\{e^{t \textbf{v}^{\perp}}(0,0) |: \|e^{t \textbf{v}^{\perp}}(0,0)\| \leq 1\right\}.$$
We want to give an upper bound on the length $\gamma$. Repeating the argument from the proof above, we can again define
a potential $A: 2 \mathbb{D} \rightarrow \mathbb{R}$ via $\nabla A = \textbf{v}^{\perp}$. Since $\textbf{v}^{\perp}$ is the gradient
flow of $A$, we see that
$$ \sup_{x \in \gamma}A(x) - \inf_{x \in \gamma}A(x) \geq c_1 |\gamma|.$$

\begin{figure}[h!]
\begin{center}
\begin{tikzpicture}[rotate=90, scale = 1]
\draw[ thick] (0,0) to[out=100,in=260] (-0.5,2.05);  
\draw (0,0) circle (60pt);
\draw[ thick] (0,0) to[out=280,in=90] (0.2,-2.1);  
\filldraw (0,0) circle (2pt);
\node at (-0.4,-0.3) {$(0,0)$};
\node at (0.2,0.2) {\Large $\gamma$};
\draw[dashed, thick] (-2,-2.5) to[out=10,in=180] (2.2,-1.5-0.5);  
\draw[dashed, thick] (-2,0) to[out=10,in=180] (2.2,1);  
\draw[dashed, thick] (-2,1.3) to[out=0,in=180] (2.2,2);  
\draw[dashed, thick] (-2,0.5) to[out=10,in=190] (2.2,1.5);  
\draw[dashed, thick] (-2,-2) to[out=0,in=190] (2.2,-1.5); 
\draw[dashed, thick] (-2,-1.5) to[out=0,in=190] (2.2,-1);  
\draw[dashed, thick] (-2,-1) to[out=0,in=190] (2.2,-0.5);   
\draw[dashed, thick] (-2,1.9) to[out=-10,in=180] (2.2,2.5);  
\end{tikzpicture}
\end{center}
\caption{Definition of $\Omega$: starting flows along $\textbf{v}$ at every point in $\gamma$ until they leave a unit disk around their starting point.}
\end{figure}
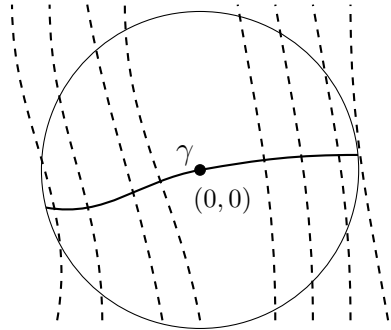
We repeat our previous argument on the domain $\Omega \subset 2\mathbb{D}$
$$ \Omega = \bigcup_{x \in \gamma}{\left\{   e^{t \textbf{v}^{}}x: \|e^{t \textbf{v}}x - x\| \leq 1\right\}}.$$
Since $A$ is constant along the flow of $\textbf{v}$, we have 
$$ \sup_{x \in \Omega}A(x) - \inf_{x \in \Omega}A(x) = \sup_{x \in \gamma}A(x) - \inf_{x \in \gamma}A(x)  \geq c_1 |\gamma|.$$ 
The coarea formula applied to $A:\Omega \rightarrow \mathbb{R}$ implies
$$ \int_{ \inf_{x \in \Omega}A(x)}^{\sup_{x \in \Omega}A(x)}{\mathcal{H}^1(A^{-1}(t))dt} =  \int_{\Omega}{| \nabla A |} = \int_{\Omega}{\| \textbf{v} \|} \leq  \int_{2\mathbb{D}}{\| \textbf{v} \|} \leq 4\pi c_2.$$
By assumption, we have $\mbox{for all} ~ t \in A(\Omega)$ 
$$\qquad \mathcal{H}^1(A^{-1}(t)) \geq 2\mathcal{L},$$
where the constant 2 is because flows both forward and backward in time are allowed.
Combining all these statements gives
$$c_1 |\gamma| 2 \mathcal{L} \leq \left(\sup_{x \in \Omega}A(x) - \inf_{x \in \Omega}A(x)\right) 2 \mathcal{L} \leq \int_{ \inf_{x \in \Omega}A(x)}^{\sup_{x \in \Omega}A(x)}{\mathcal{H}^1(A^{-1}(t))dt} \leq 4\pi c_2$$
which implies the statement
$$ |\gamma| \leq 2 \pi \frac{c_2 / c_1}{\mathcal{L}}.$$
This is the total length of $\gamma$, which leaves $\Omega$ in two points; we take the shorter path. $\qed$

\subsection{Proof of Theorem 3} 
\begin{proof}[Proof of Theorem 3] The Poincar\'{e}-Bendixson theorem coupled with the condition $\|\textbf{v}\| \geq c_1$ implies that the length of the flow until it exits the unit disk has finite length $\ell < \infty$. The same reasoning implies that if we reverse the flow of time by replacing $\textbf{v}$ with $-\textbf{v}$, the
new escape length $\ell_{-}$ is also finite. We consider $\gamma$ on the time of existence
$[-\ell_{-}, \ell]$. Furthermore, the flow has no self-intersections and $\gamma:[-\ell_{-}, \ell] \rightarrow \mathbb{R}^2$
is injective and thus splits the unit ball into two connected components $\Omega_1, \Omega_2$. There exists $i \in \left\{1, 2\right\}$ such that
$$ \mathcal{H}^1(\partial \Omega_i) \leq \ell^{-}+\ell+\pi$$
and we assume w.l.o.g. that $i=1$. 
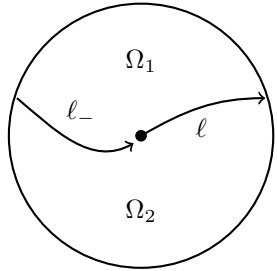
\begin{figure}[h!]
\begin{center}
\begin{tikzpicture}[scale = 1]
\filldraw (0,0) circle (2pt);
\draw [thick] (0,0) circle (50pt);
\draw[->, thick] (0,0) to[out=30,in=180] (1.65,0.5);  
\draw[->, thick] (-1.65,0.5) to[out=320,in=210] (-0.1,-0.1);  
\node at (0.8,0.1) {$\ell$};
\node at (-0.8,0.3) {$\ell_{-}$};
\node at (0,1) {$\Omega_1$};
\node at (0,-1) {$\Omega_2$};
\end{tikzpicture}
\end{center}
\caption{$\Omega_1$ is the half of the unit disk with $\Omega_1 \cap \partial \mathbb{D}$ having length at most $\pi$.}
\end{figure}

Now we consider the orthogonal vector field 
$$\textbf{v}^{\perp}(x,y) = (-v_2(x,y), v_1(x,y))$$ along $\gamma$. Trivially, $ \| \textbf{v}^{\perp}\| =\| \textbf{v}^{}\|  \geq c_1$. We split the
boundary of $\Omega_1$ into the part contained in the interior of the open unit disk $\mathbb{D}$ and the the part contained in the boundary
of the open unit disk $\mathbb{D}$, i.e.
$$ \partial \Omega_1 = (\partial \Omega_1 \cap \mathbb{D}) \cup  (\partial \Omega_1 \cap \partial \mathbb{D}).$$
On $\partial \Omega_1 \cap \mathbb{D}$ the vector field $\textbf{v}^{\perp}$ is a multiple of the normal vector
and the map
\begin{align*}
n: \partial \Omega_1 \cap \mathbb{D} &\rightarrow \mathbb{R} \\
x &\rightarrow \left\langle \nu(x), \textbf{v}^{\perp} \right\rangle
\end{align*}
is a continuous function satisfying $ |n(x)| = \|\textbf{v}^{\perp}\| = \| \textbf{v}\| \geq c_1$. Therefore
$$ \left| \int_{\partial \Omega_1 \cap \mathbb{D}}{\left\langle \textbf{v}^{\perp}, \nu(x) \right\rangle d\sigma}
\right| \geq c_1 (\ell^{-}+\ell).$$
On the other part of the boundary, we have
$$ \left| \int_{\partial \Omega_1 \cap \partial \mathbb{D}}{\left\langle \textbf{v}^{\perp}, \nu(x) \right\rangle d\sigma}
\right| \leq  \int_{\partial \Omega_1 \cap \partial \mathbb{D}}{ \|\textbf{v}^{\perp}\| d\sigma}    \leq c_2 \pi.$$
The divergence theorem implies
$$ \int_{\Omega_1}{\di \textbf{v}^{\perp} dx dy} = \int_{\partial \Omega_1}{ \left\langle \textbf{v}^{\perp}, \nu(x) \right\rangle d\sigma} = \int_{\partial \Omega_1 \cap \mathbb{D}}{\left\langle \textbf{v}^{\perp}, \nu(x) \right\rangle d\sigma} + \int_{\partial \Omega_1 \cap \partial \mathbb{D}}{\left\langle \textbf{v}^{\perp}, \nu(x) \right\rangle d\sigma} $$
and thus
\begin{align*}
 c_1(\ell^{-} + \ell) &\leq  \left| \int_{\partial \Omega_1 \cap \mathbb{D}}{\left\langle \textbf{v}^{\perp}, \nu(x) \right\rangle d\sigma} \right| \\
&\leq \left|  \int_{\Omega_1}{\di \textbf{v}^{\perp} dx dy} \right| + \left| \int_{\partial \Omega_1 \cap \partial \mathbb{D}}{\left\langle \textbf{v}^{\perp}, \nu(x) \right\rangle d\sigma} \right| \\
&\leq \left|  \int_{\Omega_1}{\di \textbf{v}^{\perp} dx dy} \right| + c_2 \pi =  \left|  \int_{\Omega_1}{\curl \textbf{v}^{} dx dy} \right| + c_2\pi\\
&\leq c_2\pi +  \int_{\mathbb{D}}{|\curl \textbf{v}^{}| dx dy}.
\end{align*}
\end{proof}

\textbf{Acknowledgment.} The author is grateful to Francesco Di Plinio for a valuable discussion and to Peter W. Jones for explaining the concept of rip currents.

\end{document}